\documentclass[reqno]{amsart}

\usepackage{latexsym,amssymb}
\usepackage{graphicx,amsmath}
\usepackage{url,amsxtra}
\usepackage{xcolor}

\def\sn{{\rm sn}}

\title{Magnetic curves corresponding to Killing magnetic fields in ${\mathbb{E}}^3$}

\author[S.~L.~Dru\c t\u a-Romaniuc]{Simona Luiza Dru\c t\u a-Romaniuc}

\address[S.L. Dru\c t\u a-Romaniuc]{\rm
'AL. I. Cuza' University of Ia\c si\\
Department of Sciences\\
Lasc\u ar Catargi Street, no.~54\\
700107 Ia\c si, Romania\\
email: {simona.druta (at) uaic.ro}
}

\author[M.~I.~Munteanu]{Marian Ioan Munteanu}

\address[M.I. Munteanu]{\rm
'Al.I.Cuza' University of Ia\c si, Faculty of Mathematics\\
Bd. Carol I, no.~11, 700506 Ia\c si, Romania\\
\url{http://www.math.uaic.ro/~munteanu}}

\address[temporary address]{\rm
Michigan State University\\
Department of Mathematics\\
Wells Hall\\
48824-1029 East Lansing\\ USA,
email: {marian.ioan.munteanu (at) gmail.com}
}

\begin{document}

\begin{abstract}
\rm We explicitly determine all magnetic curves corresponding to the Killing magnetic fields on the
3-dimensional Euclidean space.

\smallskip


\smallskip


\end{abstract}

\maketitle

\section{Introduction}

The geodesic flow on a Riemannian manifold represents the extremals of the
least action principle, namely it is determined by the motion of a certain
physical system in the manifold. It is known that the geodesic equations are
second order non-linear differential equations and they usually appear in the
form of Euler-Lagrange equations of motion. Magnetic curves generalize geodesics.
In physics, such a curve represents a trajectory of a charged particle
moving on the manifold under the action of the magnetic field.

Let $(M, g)$ be an $n$-dimensional Riemannian manifold.
A {\em magnetic field} is a closed 2-form $F$ on $M$
and the {\em Lorentz force} of a magnetic field $F$ on $(M, g)$ is an
$(1,1)$ tensor field $\Phi$ given by
\begin{equation}
\label{Lorentzforce}
g(\Phi(X), Y) = F(X, Y), \quad \forall X,Y\in \chi(M).
\end{equation}
The {\em magnetic trajectories} of $F$ are curves $\gamma$ on $M$ that
satisfy the {\em Lorentz equation} (sometimes called the {\em Newton equation})
\begin{equation}
\label{Lorentzeq}
\nabla_{\gamma^\prime}\gamma^\prime=\Phi(\gamma^\prime).
\end{equation}

Lorentz equation generalizes the equation satisfied by the geodesics of $M$, namely
$$
\nabla_{\gamma^\prime}\gamma^\prime=0.
$$
Therefore, from the point of view of the dynamical systems, a geodesic corresponds to
a trajectory of a particle without an action of a magnetic field, while a magnetic trajectory
is {\em a flowline of the dynamical system}, associated with the magnetic field.

Since the Lorentz force is skew symmetric we have
$$
\frac{d}{dt}g(\gamma^\prime,\gamma^\prime)=2g(\nabla_{\gamma^\prime}\gamma^\prime,\gamma^\prime)=0,
$$
so the magnetic curves (trajectories) have constant speed
$v(t) =||\gamma^\prime|| = v_0$. When the magnetic curve $\gamma(t)$ is
arc length parametrized $(v_0 = 1)$, it is called a {\em normal magnetic curve}.

Recall that a vector field $V$ on $M$ is {\em Killing} if and only if it satisfies
the Killing equation:
\begin{equation}
\label{Kill_eq:dm}
g(\nabla_YV,Z)+g(\nabla_ZV,Y)=0
\end{equation}
for every vector fields $Y, Z$ on $M$, where $\nabla$ is
the Levi Civita connection on $M$.

A typical example of uniform magnetic fields is obtained by multiplying the volume form on a Riemannian surface
by a scalar $s$ (usually called {\em strength}). When the surface is of constant Gaussian curvature $K$,
trajectories of such magnetic fields are well known. More precisely, on the sphere ${\mathbb{S}}^2(K)$, $K>0$,
trajectories are small (Euclidean) circles of radius $(s^2+K)^{-1/2}$, on the Euclidean plane they
are circles and the period of motion equals to $\frac{2\pi}s$, while, on a hyperbolic plane ${\mathbb{H}}^2(-K)$,
$K>0$, trajectories can be either closed curves (when $|s|>\sqrt{K}$), or open curves.
Moreover, when $|s|=\sqrt{K}$ normal trajectories are horocycles (see e.g. \cite{{Com87:dm},{Sunada}}).

This problem was extended also for different ambient spaces.
For example, if the ambient is a complex space form,
K\"ahler magnetic fields are studied (see \cite{Ada95}),
in particular, explicit trajectories for K\"ahler magnetic fields are found in
the complex projective space ${\mathbb{CP}}^n$  \cite{Ada94}.
K\"ahler magnetic fields appear in theoretical and mathematical physics, varying
from quantum field theory and string theory to general relativity.

If the ambient is a contact manifold, the fundamental 2-form defines the so-called {\em contact magnetic field}.
Interesting results are obtained when the manifold is Sasakian, namely the angle between the velocity of a normal
magnetic curve and the Reeb vector field is constant (see \cite{Cabrerizo1}). Moreover, explicit description for
normal flowlines of the contact magnetic field on a 3-dimensional Sasakian manifold is known \cite{Cabrerizo1}.

In the case of a 3-dimensional Riemannian manifold $(M,g)$, 2-forms and vector fields may be identified
via the Hodge star operator $\star$ and the volume form $dv_g$ of the manifold.
Thus, magnetic fields mean divergence free vector fields (see e.g. \cite{Cabrerizo}).
In particular, Killing vector fields define an important class of magnetic fields, called
{\em Killing magnetic fields}.
It is known that geodesics can be defined as extremal curves for the action energy functional.
A variational approach to describe Killing magnetic flows in spaces of constant curvature
is given in \cite{BarrosRomero}.

Note that, one can define on $M$ the {\em cross product}
of two vector fields $X,Y\in\chi(M)$ as follows
$$
g(X\times Y,Z)=dv_g(X,Y,Z),\quad \forall Z\in\chi(M).
$$
If $V$ is a Killing vector field on $M$, let $F_V=\iota_Vdv_g$ be
the corresponding Killing magnetic field. By $\iota$ we denote the inner product.
Then, the Lorentz force of $F_V$ is
(see \cite{Cabrerizo})
$$
\Phi(X)=V\times X.
$$
Consequently, the Lorentz force equation (\ref{Lorentzeq}) can be
written as
$$
\nabla_{\gamma^\prime}\gamma^\prime=V\times \gamma^\prime.
$$

In what follows we consider the 3-dimensional Euclidian space ${\mathbb{E}}^3$,
endowed with the usual scalar product $\langle~,~\rangle$ .

The fundamental solutions of \eqref{Kill_eq:dm} are $\{\partial_x, \partial_y,
\partial_z, -y \partial_x + x \partial_y, -z \partial_y + y \partial_z, z \partial_x - x \partial_z\}$ and
they give a basis of Killing vector fields on ${\mathbb{E}}^3$. Here $x, y, z$ denote the global
coordinates on ${\mathbb{E}}^3$ and ${\mathbb{R}}^3 = {\rm span}\{\partial_x, \partial_y,
\partial_z\}$ is regarded as a vector space.

The easiest example is to consider the Killing vector field $\xi_0=\partial_z$.
(Similar discussions can be made for $\partial_x$ and $\partial_y$, respectively.) Its trajectories
are helices with axis $\partial_z$, namely $t\mapsto(x_0+a\cos t,y_0+a \sin t, z_0+bt)$, where $(x_0,y_0,z_0)\in{\mathbb{R}}^3$
and $a,b\in{\mathbb{R}}$. An interesting fact is that Lancret curves (i.e. general helices) in ${\mathbb{E}}^3$
are characterized by the following property (in our framework): they are magnetic trajectories
associated with magnetic fields parallel to their axis. A similar result, relating Killing
magnetic fields and Lancret curves is provided on the 3-sphere (see e.g. \cite{BarrosRomero}).
Theorems of Lancret for general helices in 3-dimensional real space forms are presented in \cite{Bar97}.

In this paper we consider the following magnetic field $F_V=-(xdx+ydy)\wedge dz$ in ${\mathbb{E}}^3$,
determined by the Killing vector field
$V=-y\partial_x + x\partial_y$. The other two rotational vector fields $-z\partial_y+y\partial_z$
and $z\partial_x-x\partial_z$ give rise to analogue classifications for corresponding magnetic trajectories.
The aim of this note is to find all magnetic curves corresponding to $F_V$.
The main result we obtain is the following:

{\bf Theorem.}
{\em The magnetic trajectories of the Killing magnetic field $F_V$ are:
{\rm (a)} planar curves situated in a vertical strip; {\rm (b)} circular helices
and {\rm (c)} curve parametrized by }
$$x(t)=\rho(t)\cos\phi(t),\ y(t)=\rho(t)\sin\phi(t),\ z(t)=-\frac12\int\limits^t\rho^2(\zeta)d\zeta$$
{\em where $\rho$ and $\phi$ satisfy}
$$
\left(\frac{d\rho^2}{dt}\right)^2+P\big(\rho^2(t)\big)=0,\quad \rho^2(t)\phi'(t)={\rm constant}
$$
{\em and $P$ is a polynomial of degree $3$.}

We are able to obtain explicit solutions in case (c) and we represent some examples by using numerical
approximations for some integrals.

Recall, for later use, some basic facts on {\em normal elliptic integral of the first kind}
 (see for example \cite{BF71}):
$$
\int\limits_0^y\frac{dt}{\sqrt{(1-t^2)(1-k^2t^2)}}=\int\limits_0^\varphi\frac{d\vartheta}{\sqrt{1-k^2\sin^2\vartheta}}
=u=\sn^{-1}(y,k)=F(\varphi,k),
$$
where $y=\sin\varphi$ and $\varphi={\rm am\ } u$. The angle $\varphi$ is called {\em Jacobi amplitude} and
the function $\sn$ in known as {\em Jacobi elliptic sine}. The number $k$ is called {\em modulus} and
for applications to engineering and physics it belongs to $(0,1)$.

\section{Rotational magnetic trajectories in ${\mathbb{E}}^3$}

Let us consider the Killing vector field
$V=-y \partial_x + x\partial_y$ on ${\mathbb{E}}_3\setminus Oz$, which defines the magnetic field
$F_V=-(xdx+ydy)\wedge dz$. The Lorentz force $\Phi_V$ acts on the vector space ${\mathbb{R}}^3$
as follows:
$$
\Phi_V\partial_x=-x\partial_z,\  \Phi_V\partial_y=-y\partial_z,\
\Phi_V\partial_z=x\partial_x+y\partial_y.
$$
For the Euclidian space ${\mathbb{E}}^3$ the Lorentz force equation becomes
\begin{equation}
\label{magn}
\gamma''=V\times \gamma^\prime
\end{equation}
where the curve $\gamma:I=[0,l]\longrightarrow{\mathbb{E}}^3$, $\gamma(t)=(x(t),y(t),z(t))$ is parametrized by arc
length, namely
\begin{equation}\label{arc}
x'(t)^2+y'(t)^2+z'(t)^2=1, \quad \forall t\in I
\end{equation}
and at the moment $t=0$ it passes through the point
$(x_0,y_0,z_0)$, with the velocity $(u_0,v_0,w_0)$, such that
$$
u_0^2+v_0^2+w_0^2=1.
$$

{\bf Proof of the Theorem.}
Our aim is to determine the magnetic curves of $F_V$. The equation (\ref{magn}) yields
the following ordinary differential equations system
\begin{equation}
\label{syst}
\begin{cases}
x''=xz'\\
y''=yz'\\
z''=-(xx'+yy').
\end{cases}
\end{equation}
In order to solve it, note that from the first two equations we get a prime integral
\begin{equation}
\label{xy'}
x'y-y'x=u_0y_0-x_0v_0
\end{equation}
while from the third equation we obtain
\begin{equation}
\label{zp}
z'=-\frac{1}{2}(x^2+y^2)+\frac{1}{2}(x_0^2+y_0^2)+w_0.
\end{equation}
Notice that $z'$ cannot vanish identically (on a subinterval of $I$). Indeed, if $z'=0$ then $x'=u_0$,
$y'=v_0$ and $z=z_0$ with $u_0^2+v_0^2=1$. Hence, $x(t)=x_0+u_0t$, $y(t)=y_0+v_0t$ and combining with
\eqref{zp} we get a contradiction. It follows that one cannot have horizontal magnetic curves corresponding to $V$.

\medskip

In the sequel it is more convenient to consider cylindrical
coordinates $\{\rho,\phi,z\}$ on ${\mathbb{E}}^3\setminus Oz$. Thus, for our curve we have
$$
\begin{cases}
x=\rho(t)\cos\phi(t)\\
y=\rho(t)\sin\phi(t)\\
z=z(t)
\end{cases}
$$
where $\rho^2(t)=x^2(t)+y^2(t)$, $\rho(t)\gneq0$.

\bigskip

\textbf{Case I.}
First we study the general case, when $z'$ is not constant (equivalently $\rho$ is not constant).
The relations (\ref{xy'}) and (\ref{zp}) lead to
\begin{equation}
\label{rho}
\rho^2(t)\phi'(t)=p_0
\end{equation}
\begin{equation}
\label{zpt}
z'(t)=q_0-\frac{1}{2}\rho^2(t)
\end{equation}
where we put $p_0=x_0v_0-u_0y_0$ and $q_0=\frac{1}{2}\left(x_0^2+y_0^2\right)+w_0$.

The arc length parametrization condition (\ref{arc}), together with \eqref{zpt}, becomes
\begin{equation}\label{rhoarc}
\rho'^2(t)+\rho^2(t)\phi'^2(t)+q_0^2-q_0\rho^2(t)+\frac{1}{4}\rho^4(t)=1.
\end{equation}
Multiplying \eqref{rhoarc} by $4 \rho^2(t)$, using \eqref{rho} and denoting
$\rho^2(t)$ by $f(t)>0$, for all $t\in I$, one gets
\begin{equation}\label{F}
f'^2+f^3-4q_0f^2+4(q_0^2-1)f+4p_0^2=0.
\end{equation}

We start to study the above differential equation for some
particular values of the constants $p_0$ and $q_0$.

If $p_0=0,$ i.e. $x_0v_0=y_0u_0$ it follows that the angle $\phi$
is constant, $\phi=\phi_0$, so the magnetic trajectory is a planar curve,
with
$$
x(t)=\rho(t)\cos\phi_0,\ y(t)=\rho(t)\sin\phi_0.
$$
More precisely, the curve lies in the plane $(\sin\phi_0)x-(\cos\phi_0)y=0$.
The initial conditions expressed in cylindrical coordinates, may be written as
$$
x_0=\rho_0\cos\phi_0,\ y_0=\rho_0\sin\phi_0
$$
and the condition $x_0v_0=y_0u_0$ becomes $u_0=\zeta_0\cos\phi_0$, $v_0=\zeta_0\sin\phi_0$, for a certain $\zeta_0\in{\mathbb{R}}$.
It follows that $\zeta_0^2+\big(q_0-\frac12~\rho_0^2\big)^2=1$
\begin{equation*}
-1+\frac12~\rho_0^2\leq q_0\leq 1+\frac12~\rho_0^2
\end{equation*}
Since $\rho_0\gneq0$ it follows that $q_0>-1$.

\medskip

Let us solve the equation (\ref{F}), for three particular situations arising from the initial conditions:

\begin{enumerate}

\item [(i)] If $p_0=0$ and $q_0=0$, then the equation (\ref{F}) takes the form
$$
       f'^2(t)+f(t)\big(f(t)-2\big)\big(f(t)+2\big)=0
$$
and it has solution if and only if $f(t)\leq 2$, i.e. $\rho(t)\in (0,\sqrt{2}]$, so the magnetic
curve $\gamma$ lies inside a cylinder. In fact, being a planar curve,  $\gamma$ stays in a
vertical strip centered in $Oz$ and of width $2\sqrt{2}$.

\noindent
We have $f'(t)=\pm\sqrt{f(t)\big(4-f(t)^2\big)}$ and we consider only the {\em plus} sign (the other situation
may be treated in similar way).
Supposing $\rho_0\neq \sqrt{2}$, we have that $f$ and the
integral ${\mathcal{I}}(f)=\displaystyle\int^f_{\rho_0^2} \frac{d\zeta}{\sqrt{\zeta (4-\zeta^2)}}$ are strictly increasing functions.
Thus, the equation ${\mathcal{I}}(f)=t$ has a unique solution in the interval $(\rho_0^2,2)$, namely $f={\mathcal{J}}(t)$,
where ${\mathcal{J}}$ is the inverse function of ${\mathcal{I}}$.
Consequently, $\rho(t)=\sqrt{{\mathcal{J}}(t)}$.
In fact ${\mathcal{J}}$ may be expressed in terms of the elliptic functions. More precisely,
$$
{\mathcal{J}}(t)=\frac{2~\sn^2(t+t_0,\frac1{\sqrt{2}})}{2-\sn^2(t+t_0,\frac1{\sqrt{2}})}
$$
where $t_0$ is determined by $\sn(t_0,\frac1{\sqrt{2}})=\frac{\sqrt{2}\rho_0}{\sqrt{2+\rho_0^2}}$.

\noindent
Summarizing, the magnetic curve is given by
$$
x(t)=\sqrt{{\mathcal{J}}(t)}\cos\phi_0,\
y(t)=\sqrt{{\mathcal{J}}(t)}\sin\phi_0,\
z(t)=-\frac{1}{2}\int_0^t {\mathcal{J}}(\zeta)d\zeta.
$$

\noindent
In order to draw a picture of our curve, one can use Matlab to compute the parametrization.
The idea is to calculate the integrals numerically, as Riemann sums. See Appendix.
\label{Ii}

\item [(ii)] If $p_0=0,\ q_0=1$, then the equation (\ref{F}) becomes
$f'^2(t)+f^2(t)\big(f(t)-4\big)=0,$ from which we have that $f(t)\leq 4,$
equivalently $\rho(t)\leq 2$, so the magnetic curve $\gamma$ stays inside
a cylinder of radius 2. In fact, being planar, the curve lies in a vertical strip centered on $z$-axis.
The equation can be written in the form
$$
\frac{df}{f\sqrt{4-f}}=\pm~dt.
$$
Taking the {\em plus} sign, one gets the solution
$$
f(t)=\frac{4}{\cosh^2(t-t_0)},\quad t\in (0,t_0)
$$
where
$t_0=-\frac{1}{2}\ln\frac{2-\sqrt{4-\rho_0^2}}{2+\sqrt{4-\rho_0^2}}$.
Hence
$$
\rho(t)=\frac{2}{\cosh(t-t_0)}
$$
and the magnetic curve is parametrized by
$$
x(t)=\frac{2\cos\phi_0}{\cosh(t-t_0)},\
y(t)=\frac{2\sin\phi_0}{\cosh(t-t_0)},\
z(t)=z_0+t-2 \left(\tanh(t-t_0)+\tanh t_0\right).
$$
\end{enumerate}

We draw a picture of this (planar) curve.

\begin{figure}[hbtp]
\label{fig:sol_03}
\begin{center}
  \includegraphics[width=75mm]{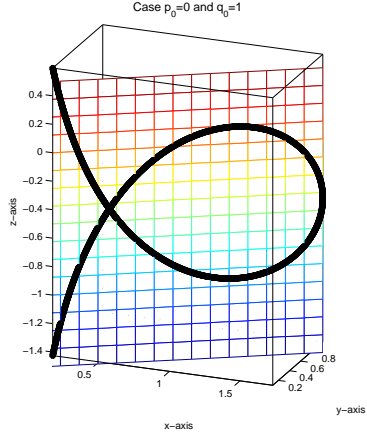}
\end{center}
\caption{$t\in{\mathbb{R}}$, $\phi_0=\frac\pi6$}
\end{figure}

Let us finalize the examination of the equation \eqref{F} for $p_0=0$. The polynomial
\linebreak
$P(f)=f\Big(f^2-4q_0f+4(q_0^2-1)\Big)$ has three
real solutions, namely $f_1=2(q_0-1)$, $f_2=2(q_0+1)$ and $f_3=0$.
If $f$ is a solution for \eqref{F}, then $P(f)$ should be negative.
Recall that $q_0>-1$. We have
\begin{enumerate}
\item [(a)] If $q_0\in(-1,1)$, then $f_1<0<f_2$.
It follows that $\rho(t)\in(0,\sqrt{2(q_0+1)})$ and the discussion is similar
as in case $q_0=0$. More precisely we have
$$
\rho(t)^2=\frac{2(1-q_0^2)~\sn^2\big(t+t_0,\sqrt{\frac{q_0+1}2}\big)}{2-(q_0+1)~\sn^2\big(t+t_0,\sqrt{\frac{q_0+1}2}\big)}
$$
where $t_0$ is defined by $\sn\big(t_0,\sqrt{\frac{q_0+1}2}\big)=\sqrt{\frac2{q_0+1}}\frac{\rho_0}{\sqrt{\rho_0^2-2(q_0-1)}}$.

\item [(b)] If $q_0>1$, then $0<f_1<f_2$.
It follows that $f(t)\in(f_1,f_2)$. Thus, the curve $\gamma$ lies between two cylinders
since $\rho(t)\in\left(\sqrt{2(q_0-1)},\sqrt{2(q_0+1)}\right)$. As before, the curve is situated
in a union of two vertical strips. Again, the discussion is similar as in case $q_0=0$.
In terms of elliptic functions, we may write
$$
\rho(t)^2=\frac{q_0^2-1}{\frac{q_0+1}2-\sn^2\big(\sqrt{\frac{q_0+1}2} t+t_0,\sqrt{\frac2{q_0+1}}\big)}
$$
where $t_0$ is defined by $\sn\big(t_0,\sqrt{\frac2{q_0+1}}\big)=\sqrt{\frac{q_0+1}2}\frac{\sqrt{\rho_0^2-2(q_0-1)}}{\rho_0}$.

\medskip
\noindent
In order to visualize an example, consider the following initial conditions:
$x_0=2$, $y_0=0$, $z_0=0$ and $u_0=0$, $v_0=0$, $w_0=1$ (this yields $p_0=0$ and $q_0=3$).

\medskip

\noindent
We will use again Matlab to compute the integrals (numerically) and to draw the picture.
\end{enumerate}

\begin{figure}[hbtp]
\label{fig:exb_pg5}
\begin{center}
  \includegraphics[width=75mm]{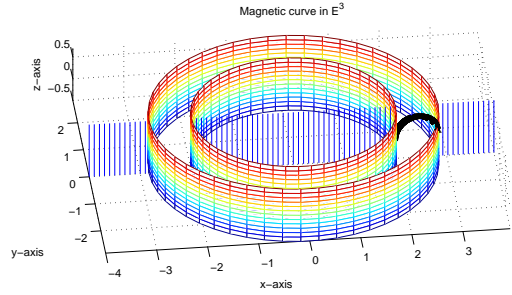}
\end{center}
\caption{$p_0=0$, $q_0=3$}
\end{figure}


Return to \eqref{F} for $p_0\neq0$ and notice that the equation
\begin{equation}
\label{eq14}
P(f)=f^3-4q_0f^2+4(q_0^2-1)f+4p_0^2=0
\end{equation}
has the discriminant
$$
\Delta=-16\big[27p_0^4+8p_0^2q_0(q_0^2-9)-16(q_0^2-1)^2\big]
$$
and the following situations appear:
\begin{itemize}
    \item the equation \eqref{eq14} has three distinct solutions iff $\Delta>0.$
    \item the polynomial $P$ has multiple roots iff $\Delta=0$.
    \item the polynomial $P$ has one real root and two complex
    conjugate roots iff $\Delta<0$.
\end{itemize}

A detailed analysis of the above situations, lead us to conclude, after
taking into account classical Vi\`ete's formulas, that the equation (\ref{F})
has solutions if and only if $\Delta>0$.

Indeed, if $\Delta<0$, let $A\in{\mathbb{C}}\setminus{\mathbb{R}}$
and $\bar A$ be the complex solutions of \eqref{eq14}, and $B$ its real
solution. Then, the ODE \eqref{F} can be rewritten as
$$
    f'(t)^2+\left(f(t)^2-2~Re(A)~t+|A|^2\right)\left(f(t)-B\right)=0,
$$
where $Re(A)$ denotes the real part of the complex number $A$.
From the third Vi\`ete's formula we conclude that $B$ should be negative,
and consequently, the previous equality cannot occur.

On the other hand, if $\Delta=0$, analyzing the coefficients one cannot have
a triple root (since $16q_0^2\neq12(q_0^2-1)$.
Hence, let $A\in{\mathbb{R}}$ be the double root, and let $B\in{\mathbb{R}}$ be
the third one. With a similar argument as above, $B$ is negative and the ODE \eqref{F} becomes
$$
    f'(t)^2+\left(f(t)-A\right)^2\left(f(t)-B\right)=0.
$$
Again, this equality cannot hold.

It follows that $\Delta$ should be (strictly) positive.
Let $A$, $B$, $C\in {\mathbb{R}}$ be the three distinct solutions of \eqref{eq14}.
The third Vi\`ete's formula yields $ABC=-4p_0^2<0$, and hence
\begin{itemize}
\item [a)] either $A,B,C$ are all negative,
\item [b)] or two of them, $A$ and $B$, are positive and the third one, $C$, is negative.
\end{itemize}
In case a) the ODE \eqref{F}: $f'(t)^2+(f(t)-A)(f(t)-B)(f(t)-C)=0$ has no solution.
This happens if and only if $q_0<-1$ and $p_0\neq0$ (for the proof use the second and the third Vi\`ete's formulas)
together with $\Delta>0$ (for example if $q_0=-3$ and $p_0=1$).

In case b), equivalently to $\Delta>0$, $q_0>-1$ and $p_0\neq 0$, the equation \eqref{F} has a solution in the
interval defined by the positive solutions $A<B$ of \eqref{eq14}. Since the function
\linebreak
${\mathcal{I}}(f)=\displaystyle\int_{A}^f\frac{d\zeta}{\sqrt{(\zeta-A)(B-\zeta)(\zeta-C)}}$ is
strictly increasing, ${\mathcal{I}}(f)=t$ has a unique solution $f$, denoted by
${\mathcal{J}}(t)$. Thus we have $\rho=\sqrt{{\mathcal{J}}(t)}$, and
$\phi(t)=\phi_0+p_0\displaystyle\int_0^t\frac{d\zeta}{{\mathcal{J}}(\zeta)}$.
In this case, the magnetic curve $\gamma$ is given by
$$
x(t)=\sqrt{{\mathcal{J}}(t)}\cos\Big(\phi_0+p_0\int_0^t\frac{d\zeta}{{\mathcal{J}}(\zeta)}\Big),\
y(t)=\sqrt{{\mathcal{J}}(t)}\sin\Big(\phi_0+p_0\int_0^t\frac{d\zeta}{{\mathcal{J}}(\zeta)}\Big),
$$
$$
z(t)=z_0+q_0t-\frac{1}{2}\int_0^t{\mathcal{J}}(\zeta)d\zeta.
$$
We may express $\rho$ in terms of elliptic functions, namely
$$
\rho(t)^2=\frac{Ak^2-C\sn^2(rt+t_0,\frac1k)}{k^2-\sn^2(rt+t_0,\frac1k)}
$$
where $k^2=\frac{B-C}{B-A}$, $r=\frac{\sqrt{B-C}}2$, and $\sn(t_0,\frac1k)=k\sqrt{\frac{\rho_0^2-A}{\rho_0^2-C}}$.

In the Appendix we will draw a picture (using the same technique in Matlab as before) corresponding to the following data: $p_0=\frac{\sqrt{2\sqrt{6}-3}}2$,
$q_0=\frac{3-\sqrt{6}} 2$, for which we have $A=1$, $B=2$ and $C=3-2\sqrt{6}$.
\label{ex_ABC}

The situation $I(f)=-t$ can be treated in similar way.

Finally, notice that for $q_0=-1$ we get $\Delta=-16p_0^2(27p_0^2+64)$ and this case was discussed above.
\vskip2mm

\textbf{Case II.} Now, let us study the remaining case when $z'(t)=w_0\neq0$.
We immediately have that
$$
z(t)=z_0+tw_0,
$$
and from (\ref{zp}) we obtain
\begin{equation}\label{xy}
x^2+y^2=x_0^2+y_0^2.
\end{equation}
This means that the magnetic trajectory $\gamma$ lies on the circular cylinder of radius
\linebreak
$\rho_0=\sqrt{x_0^2+y_0^2}$.

Two subcases must be discussed: $w_0<0$ and $w_0>0$.

\begin{description}
\item [II.1] In the case when $w_0<0$ the magnetic curve is given by
\begin{equation}
\label{1}
\begin{cases}
x(t)=x_0\cos(\sqrt{-w_0}~t)+\frac{u_0}{\sqrt{-w_0}}\sin(\sqrt{-w_0}~t)\\[2mm]
y(t)=y_0\cos(\sqrt{-w_0}~t)+\frac{v_0}{\sqrt{-w_0}}\sin(\sqrt{-w_0}~t)\\[2mm]
z(t)=z_0+tw_0.
\end{cases}
\end{equation}
This curve is a helix around the above cylinder.

\noindent
At this point, we have to find which are the initial conditions leading this situation.
To do this, using \eqref{xy} and \eqref{1}, we should have the following relations
\begin{equation*}
w_0=-\frac2{\rho_0^2+\sqrt{\rho_0^4+4}}\ ,\ u_0=\varepsilon\rho_0\sqrt{-w_0}\sin\phi_0,\
    v_0=-\varepsilon\rho_0\sqrt{-w_0}\cos\phi_0
\end{equation*}
where $\rho_0$ and $\phi_0$ have the usual meaning and $\varepsilon=\pm1$.

\item [II.2] If $w_0>0$, the ODE system \eqref{syst} has the following solution
\begin{equation*}
\begin{cases}
x(t)=x_0\cosh(\sqrt{w_0}~t)+\frac{u_0}{\sqrt{w_0}}\sinh(\sqrt{w_0}~t)\\
y(t)=y_0\cosh(\sqrt{w_0}~t)+\frac{v_0}{\sqrt{w_0}}\sinh(\sqrt{w_0}~t)\\
z(t)=z_0+tw_0
\end{cases}
\end{equation*}
but in this case the condition \eqref{xy} is satisfied if and only
if $x_0=y_0=0$ and $u_0=v_0=0$. This situation cannot occur.
\end{description}

\section{Review on the classical magnetic field on ${\mathbb{E}}^3$}

As we have already said in Introduction, the best known example of magnetic fields
in the Euclidean space ${\mathbb{E}}^3$ is furnished by the 2-form $F_0=dx\wedge dy$,
corresponding to the Killing vector field $\xi_0=\frac\partial{\partial z}$.

In this section we consider the Killing magnetic field $F_\xi=s~F_0=s~dx\wedge dy$,
determined by the Killing vector field $\xi=s~\xi_0=s~\partial_z$ on ${\mathbb{E}}^3$,
where $s\neq0$ is an arbitrary constant.
We briefly describe its magnetic curves.

The action of the Lorentz force $\Phi_\xi$ on the vector space ${\mathbb{R}}^3$ is given by:
$$
\Phi_\xi\partial_x=s\partial_y,\ \Phi_\xi\partial_y=-s\partial_x,\
\Phi_\xi\partial_z=0.
$$

Solving the Lorentz force equation $\gamma_s''=\Phi_\xi(\gamma_s')$,
we obtain the family of magnetic curves $\gamma_s(t)=(x(t),y(t),z(t))$, parametrized by
$$
\begin{cases}
x(t)=\frac{u_0}{s}\sin(st)+\frac{v_0}{s}\cos(st)+x_0-\frac{v_0}{s}\\
y(t)=-\frac{u_0}{s}\cos(st)+\frac{v_0}{s}\sin(st)+y_0+\frac{u_0}{s}\\
z(t)=w_0t+z_0.
\end{cases}
$$

Write the first Fr\'enet equation
$$
\gamma''=\kappa N
$$
where $\kappa$ is the curvature and $N$ is the normal of the curve.
Using the equation \eqref{magn} we obtain that the square of the curvature is
$$
\kappa^2=s^2(1-w_0^2).
$$
Moreover, classical computations give the torsion $\tau=sw_0.$

Notice that even both the curvature $\kappa$ and the torsion $\tau$ depend on the
strength $s$, the ratio $\tau\over\kappa$ does not.

We conclude with some comments:
\begin{enumerate}
\item [i)] If $w_0=0$ the curvature is $\kappa=s$ and the torsion is $\tau=0,$
so the magnetic line is a (planar) circle.

\item [ii)] If $w_0=\pm 1,$ then $\kappa=0,\ \tau=\pm s$, so the magnetic
curves are vertical lines.

\item [iii)] In other cases the magnetic curves are circular helices.
\end{enumerate}

\section{Appendix}
In this section we present a Matlab program in order to compute,
by numerical approximation of the involved integrals, the parametrization of magnetic curve
obtained in case I (i) from page \pageref{Ii}. Since the curve is planar we consider $\phi_0=0$.

{
\begin{verbatim}
clear all
%%% Compute the integral I(f) as a Riemann sum

rho0=1.41;
f_max=2;
N=1000;
L=(f_max-rho0^2)/N;

for K=1:N+1
    a=0.001;
    b=rho0^2+(K-1)*L;
    n=1000;
    h=(b-a)/n;

    k=0:n-1;
    x=a+k*h;
    f=1./sqrt(x.*(4-x.^2));

    I(K)=h*sum(f);
    J(K)=b;
end

%% \phi_0=0

xx=sqrt(J);
%yy=0*J;

zz(1)=0;
for K=1:N
    zz(K+1)=zz(K)-0.5*(I(K+1)-I(K))*J(K);
end

%% the curve is planar
plot(xx,zz,'g-')
text(0.25,-0.75,'\rho_0=1.41','Color','g')
hold on
\end{verbatim}
}

Representation of the magnetic curves depending on the initial position:

\pagebreak

\begin{figure}[hbtp]
\label{fig:Ii}
\begin{center}
  \includegraphics[width=85mm]{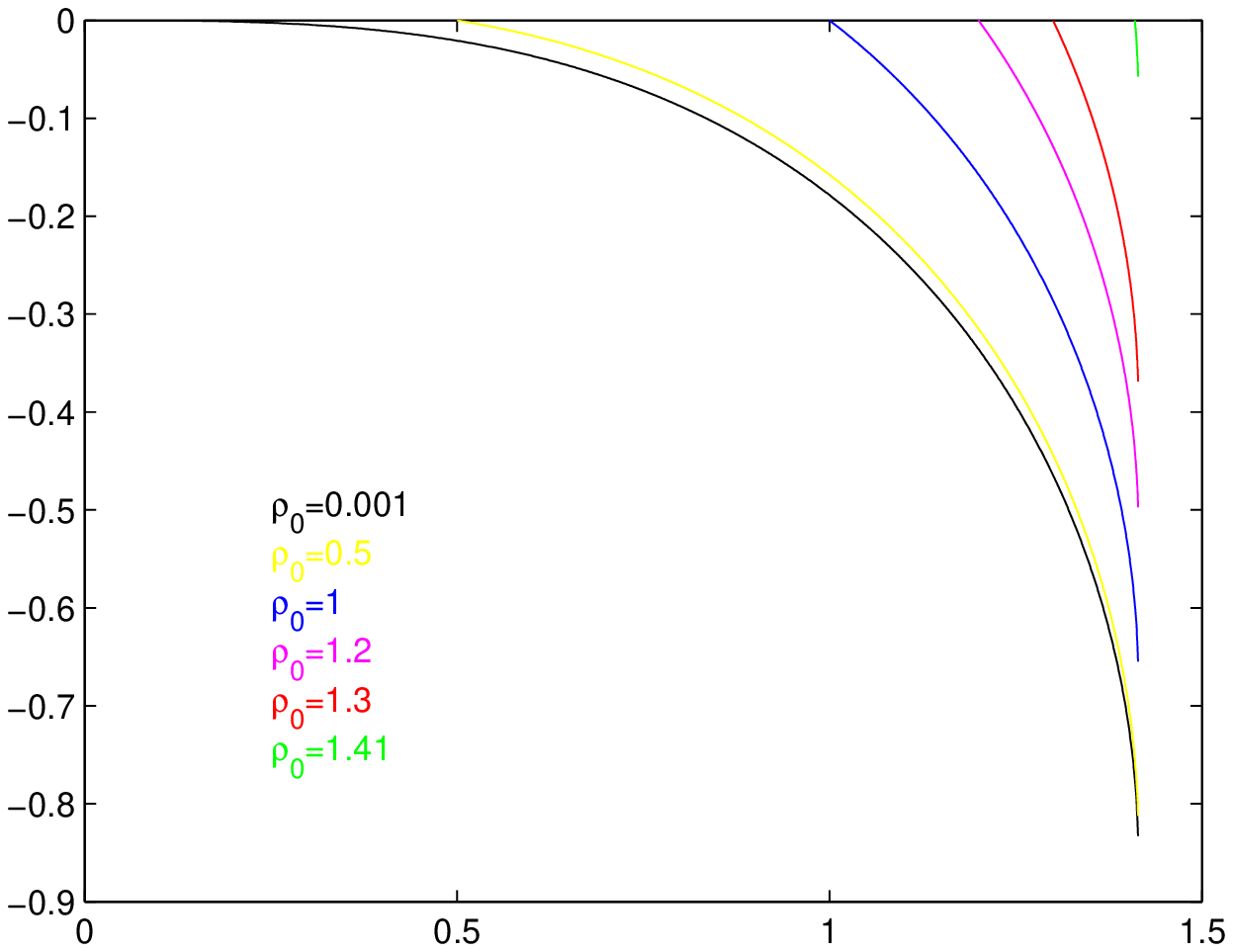}
\end{center}
\caption{I (i)}
\end{figure}

Using the previous Matlab program adapted to the example furnished at page~\pageref{ex_ABC},
and for the initial data $\phi_0=0$ and $z_0=0$, we can represent the corresponding magnetic curves:

\begin{figure}[hbtp]
\label{fig:Ii}
\begin{center}
  \includegraphics[width=85mm]{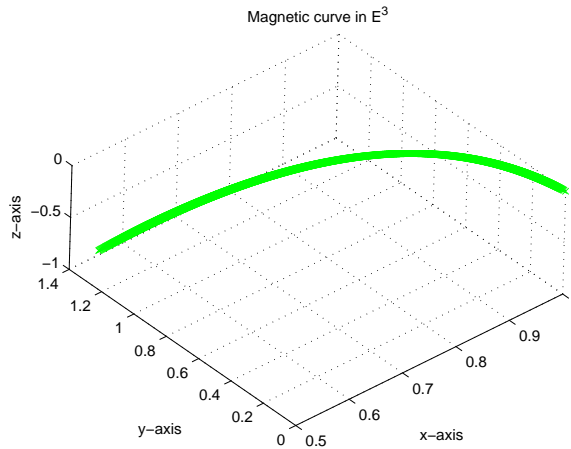}
\end{center}
\caption{$x_0=1$, $y_0=0$, $z_0=0$, $u_0=0$, $v_0=\sqrt{\frac{2\sqrt{6}-3}2}$ and $w_0=\frac{2-\sqrt{6}}{2}$}
\end{figure}

{\bf Acknowledgements.}
The first authors is a postdoctoral researcher in the framework of the program POSDRU 89/1.5/S/49944,
'AL. I. Cuza' University of Ia\c si, Romania. The second author is supported by a Fulbright Grant no. 498 at the Michigan State University, USA.

\bibliographystyle{model1-num-names}

\end{document}